\documentclass[11pt,a4paper]{amsart}

\usepackage{calc,amsthm,amsfonts,amsmath,graphicx}
\usepackage[numbers,sort&compress]{natbib}
\usepackage[lmargin=31mm,rmargin=31mm,tmargin=30mm,bmargin=30mm,footskip=10mm]{geometry}

\theoremstyle{plain}
\newtheorem{theorem}{Theorem}
\newtheorem{lemma}{Lemma}

\theoremstyle{definition}

\newcommand{\thmlabel}[1]{\label{thm:#1}}
\newcommand{\thmref}[1]{Theorem~\ref{thm:#1}}

\newcommand{\lemlabel}[1]{\label{lem:#1}}
\newcommand{\lemref}[1]{Lemma~\ref{lem:#1}}

\newcommand{\figlabel}[1]{\label{fig:#1}}
\newcommand{\figref}[1]{Figure~\ref{fig:#1}}

\newcommand{\seclabel}[1]{\label{sec:#1}}
\newcommand{\secref}[1]{Section~\ref{sec:#1}}

\allowdisplaybreaks

\newcommand{\Figure}[4][htb]{
\begin{figure}[#1]
	\vspace*{1ex}
	\begin{center}#3\end{center}
	\vspace*{-1ex}
	\caption{\figlabel{#2}#4}
\end{figure}
}

\newcommand{\ceil}[1]{\ensuremath{\protect\lceil#1\rceil}}

\newcommand{\half}{\ensuremath{\protect\tfrac{1}{2}}}

\newcommand{\eighth}{\ensuremath{\protect\tfrac{1}{8}}}

\newcommand{\quarter}{\ensuremath{\protect\tfrac{1}{4}}}

\newcommand{\twothirds}{\ensuremath{\protect\tfrac{2}{3}}}

\newcommand{\Oh}[1]{\ensuremath{\protect\mathcal{O}(#1)}}

\newcommand{\tpw}[1]{\ensuremath{\textup{\textsf{tpw}}(#1)}}
\newcommand{\dtw}[1]{\ensuremath{\textup{\textsf{dtw}}(#1)}}
\newcommand{\tw}[1]{\ensuremath{\textup{\textsf{tw}}(#1)}}

\renewcommand{\baselinestretch}{1.135}
\begin{document}

\title{On Tree-Partition-Width}

\author{David R. Wood}

\address{\newline Departament de Matem{\`a}tica Aplicada II\newline Universitat Polit{\`e}cnica de Catalunya\newline Barcelona, Spain}

\thanks{Research supported by a Marie Curie Fellowship of the European Community under contract MEIF-CT-2006-023865, and by the projects MEC MTM2006-01267 and DURSI 2005SGR00692.}

\email{david.wood@upc.es}

\date{\today}

\subjclass{05C70 (graph theory: factorization, matching, covering and packing)}

\begin{abstract}
A \emph{tree-partition} of a graph $G$ is a proper partition of its vertex set into  `bags', such that identifying the vertices in each bag produces a forest. The \emph{tree-partition-width} of $G$ is the  minimum number of vertices in a bag in a tree-partition of $G$. An anonymous referee of the paper by Ding and Oporowski [\emph{J.\ Graph Theory}, 1995] proved that every graph with tree-width $k\geq3$ and maximum degree $\Delta\geq1$ has tree-partition-width at most $24k\Delta$. We prove that this bound is within a constant factor of optimal. In particular, for all $k\geq3$ and for all sufficiently large $\Delta$, we construct a graph with tree-width $k$, maximum degree $\Delta$, and tree-partition-width at least $(\eighth-\epsilon)k\Delta$. Moreover, we slightly improve the upper bound to $\frac{5}{2}(k+1)(\frac{7}{2}\Delta-1)$ without the restriction that $k\geq3$. 
\end{abstract}

\keywords{graph, tree-partition, tree-partition-width, tree-width}

\maketitle

\section{Introduction}


A graph\footnote{All graphs considered are undirected, simple, and finite. Let $V(G)$ and $E(G)$ respectively be the vertex set and edge set of a graph $G$. Let $\Delta(G)$ be the maximum degree of $G$.} $H$ is a \emph{partition} of a graph $G$ if:
\begin{itemize}
\item each vertex of $H$ is a set of vertices of $G$ (called a \emph{bag}), 
\item every vertex of $G$ is in exactly one bag of $H$, and 
\item distinct bags $A$ and $B$ are adjacent in $H$ if and only if some edge of $G$ has one endpoint in $A$ and the other endpoint in $B$. 
\end{itemize}
The \emph{width} of a partition is the maximum number of vertices in a bag. Informally speaking, the graph $H$ is obtained from a proper partition of $V(G)$ by identifying the vertices in each part, deleting loops, and replacing parallel edges by a single edge. 


If a forest $T$ is a partition of a graph $G$, then $T$ is a \emph{tree-partition} of $G$. The \emph{tree-partition-width}\footnote{Tree-partition-width has also been called \emph{strong tree-width} \cite{Seese85, BodEng-JAlg97}.} of $G$, denoted by  \tpw{G}, is the minimum width of a tree-partition of $G$. Tree-partitions were independently introduced by \citet{Seese85} and \citet{Halin91}, and have since been widely investigated \citep{BodEng-JAlg97, Bodlaender-DMTCS99, DO-JGT95, DO-DM96, Edenbrandt86, Wood-JGT06}. Applications of tree-partitions include graph drawing \citep{GLM-CGTA05, DSW-CGTA, WoodTelle-NYJM07, DMW-SJC05}, graph colouring \citep{WoodBarat-arXiv}, partitioning graphs into subgraphs with only small components \citep{ADOV-JCTB03}, monadic second-order logic \citep{KuskeLohrey05}, and network emulations \citep{BvL-IC86, FF82, Bodlaender-IC90, Bodlaender-IPL88}. Planar-partitions and other more general structures have also recently been studied \citep{DiestelKuhn-DAM05,WoodTelle-NYJM07}.


What bounds can be proved on the tree-partition-width of a graph? Let \tw{G} denote the tree-width\footnote{A graph is \emph{chordal} if every induced cycle is a triangle. The \emph{tree-width} of a graph $G$ can be defined to be the minimum integer $k$ such that $G$ is a subgraph of a chordal graph with no clique on $k+2$ vertices. This parameter is particularly important in algorithmic and structural graph theory; see \citep{Bodlaender-TCS98, Reed-AlgoTreeWidth03} for surveys.} of a graph $G$. \citet{Seese85} proved the lower bound,
$$2\,\tpw{G}\geq\tw{G}+1.$$
In general, tree-partition-width is not bounded from above by any function solely of tree-width. For example, wheel graphs have bounded tree-width and unbounded tree-partition-width \citep{BodEng-JAlg97}. However, tree-partition-width is bounded for graphs of bounded tree-width \emph{and} bounded degree \cite{DO-JGT95, DO-DM96}. The best known upper bound is due to an anonymous referee of the paper by \citet{DO-JGT95}, who proved that $$\tpw{G}\leq24\,\tw{G}\,\Delta(G)$$ whenever $\tw{G}\geq3$ and $\Delta(G)\geq1$. Using a similar proof, we make the following improvement to this bound without the restriction that $\tw{G}\geq3$.

\begin{theorem}
\thmlabel{Tweaking}
Every graph $G$ with tree-width $\tw{G}\geq1$ and maximum degree $\Delta(G)\geq1$ has tree-partition-width  $$\tpw{G}<\tfrac{5}{2}\big(\tw{G}+1\big)\big(\tfrac{7}{2}\,\Delta(G)-1\big).$$
\end{theorem}

\thmref{Tweaking} is proved in \secref{UpperBound}. Note that  \thmref{Tweaking} can be improved in the case of chordal graphs. In particular, a simple extension of a result by \citet{DMW-SJC05} implies that $$\tpw{G}\leq\tw{G}\big(\Delta(G)-1\big)$$ for 
every chordal graph $G$ with $\Delta(G)\geq2$; see \citep{Wood-JGT06} for a simple proof. Nevertheless, the following theorem proves that \Oh{\tw{G}\,\Delta(G)} is the best possible upper bound, even for chordal graphs.

\begin{theorem}
\thmlabel{GeneralLowerBound}
For every $\epsilon>0$ and integer $k\geq3$, for every sufficiently large integer $\Delta\geq\Delta(k,\epsilon)$, for infinitely many values of $N$, there is a chordal graph $G$ with $N$ vertices, tree-width $\tw{G}\leq k$, maximum degree $\Delta(G)\leq\Delta$, and tree-partition-width $$\tpw{G}\geq(\eighth-\epsilon)\,\tw{G}\,\Delta(G).$$
\end{theorem}

\thmref{GeneralLowerBound} is proved in \secref{GeneralLowerBound}. Note that \thmref{GeneralLowerBound} is for $k\geq3$. For $k=1$, every tree is a tree-partition of itself with width $1$. For $k=2$, we prove that the upper bound \Oh{\Delta(G)} is again best possible; see \secref{SeriesParallel}.

\section{Upper Bound}
\seclabel{UpperBound}

In this section we prove \thmref{Tweaking}. The proof relies on the following separator lemma by \citet{RS-GraphMinorsII-JAlg86}.

\begin{lemma}[\citep{RS-GraphMinorsII-JAlg86}]
\lemlabel{Separator}
For every graph $G$ with tree-width at most $k$, for every set $S\subseteq V(G)$, there are edge-disjoint subgraphs $G_1$ and $G_2$ of $G$ such that $G_1\cup G_2=G$, $|V(G_1)\cap V(G_2)|\leq k+1$, and $|S-V(G_i)|\leq\frac{2}{3}|S-(V(G_1)\cap V(G_2))|$ for each $i\in\{1,2\}$.
\end{lemma}

\thmref{Tweaking} is a corollary of the following stronger result. 

\begin{lemma}
\lemlabel{Tweaking}
Let $\alpha:=1+1/\sqrt{2}$ and $\gamma:=1+\sqrt{2}$.
Let $G$ be a graph with tree-width at most $k\geq1$ and maximum degree at most $\Delta\geq1$. Then $G$ has tree-partition-width  $$\tpw{G}\leq\gamma(k+1)(3\gamma\Delta-1)\enspace.$$ 
Moreover, for each set $S\subseteq V(G)$ such that 
$$(\gamma+1)(k+1)\leq|S|\leq3(\gamma+1)(k+1)\Delta,$$ 
there is a tree-partition of $G$ with width at most 
$$\gamma(k+1)(3\gamma\Delta-1),$$ 
such that $S$ is contained in a single bag containing at most $\alpha|S|-\gamma(k+1)$ vertices. 
\end{lemma}

\begin{proof}
We proceed by induction on $|V(G)|$. 

\textbf{Case 1.}~$|V(G)|<(\gamma+1)(k+1)$: Then no set $S$ is specified, and the tree-partition in which all the vertices are in a single bag satisfies the lemma. Now assume that $|V(G)|\geq(\gamma+1)(k+1)$, and without loss of generality, $S$ is specified. 

\textbf{Case 2.}~$|V(G)-S|<(\gamma+1)(k+1)$: Then the tree-partition in which $S$ is one bag and $V(G)-S$ is another bag satisfies the lemma. Now assume that $|V(G)-S|\geq(\gamma+1)(k+1)$. 

\textbf{Case 3.}~$|S|\leq3(\gamma+1)(k+1)$: Let $N$ be the set of vertices in $G$ that are adjacent to some vertex in $S$ but are not in $S$. Then $|N|\leq\Delta|S|\leq3(\gamma+1)(k+1)\Delta$. If $|N|<(\gamma+1)(k+1)$ then add arbitrary vertices from $V(G)-(S\cup N)$ to $N$ until $|N|\geq(\gamma+1)(k+1)$. This is possible since $|V(G)-S|\geq(\gamma+1)(k+1)$. 

By induction, there is a tree-partition of $G-S$ with width at most $\gamma(k+1)(3\gamma\Delta-1)$, such that $N$ is contained in a single bag.
Create a new bag only containing $S$. Since all the neighbours of $S$ are in a single bag, we obtain a tree-partition of $G$. ($S$ corresponds to a leaf in the pattern.)\ 
Since $|S|\geq(\gamma+1)(k+1)$, 
it follows that $|S|\leq\alpha|S|-\gamma(k+1)$ as desired. 
Now $|S|\leq3(\gamma+1)(k+1)<\gamma(k+1)(3\gamma\Delta-1)$. 
Since the other bags do not change we have the desired tree-partition of $G$.

\textbf{Case 4.}~$|S|\geq3(\gamma+1)(k+1)$: By \lemref{Separator}, there are edge-disjoint subgraphs $G_1$ and $G_2$ of $G$ such that $G_1\cup G_2=G$, $|V(G_1)\cap V(G_2)|\leq k+1$, and $|S-V(G_i)|\leq\frac{2}{3}|S-(V(G_1)\cap V(G_2))|$ for each $i\in\{1,2\}$. Let $Y:=V(G_1)\cap V(G_2)$. Let $a:=|S\cap Y|$ and $b:=|Y-S|$. Thus $a+b\leq k+1$. Let $p_i:=|(S\cap V(G_i))-Y|$. Then $p_1\leq 2p_2$ and $p_2\leq 2p_1$. Let $S_i:=(S\cap V(G_i))\cup Y$. Note that $|S_i|=p_i+a+b$. 

\begin{figure}[!ht]
\begin{center}
\includegraphics{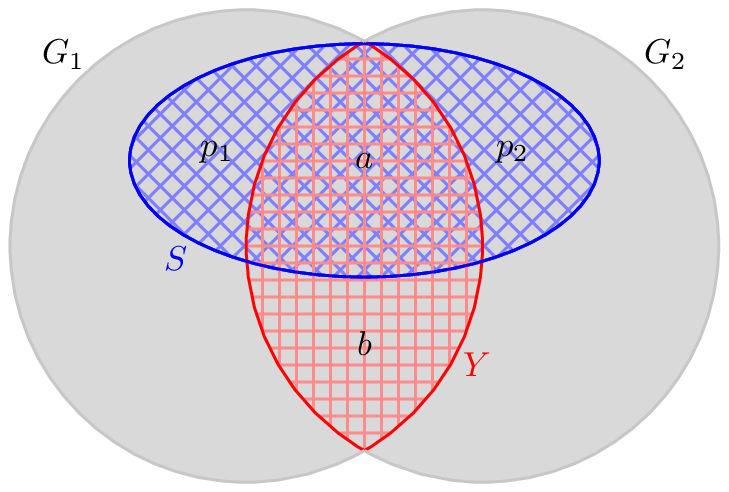}
\caption{Illustration of Case 4.}
\end{center}
\end{figure}

Now $p_1+p_2+a=|S|\geq3(\gamma+1)(k+1)$. Thus 
$3p_i+a\geq3(\gamma+1)(k+1)$
and
$3p_i+3a+3b\geq3(\gamma+1)(k+1)$.
That is,
$|S_i|\geq(\gamma+1)(k+1)$ for each $i\in\{1,2\}$.

Now $p_1+p_2+a\leq3(\gamma+1)(k+1)\Delta$. Thus 
$\frac{3}{2}p_i+a\leq3(\gamma+1)(k+1)\Delta$
and
$p_i\leq2(\gamma+1)(k+1)\Delta$. Thus
$p_i+a+b\leq2(\gamma+1)(k+1)\Delta+(k+1)$.
Hence $|S_i|=p_i+a+b<3(\gamma+1)(k+1)\Delta$. 

Thus we can apply induction to the set $S_i$ in the graph $G_i$ for each $i\in\{1,2\}$. We obtain a tree-partition of $G_i$ with width at most $\gamma(k+1)(3\gamma\Delta-1)$, such that $S_i$ is contained in a single bag $T_i$ containing at most $\alpha|S_i|-\gamma(k+1)$ vertices. 

Construct a partition of $G$ by uniting $T_1$ and $T_2$. Each vertex of $G$ is in exactly one bag since $V(G_1)\cap V(G_2)=Y\subseteq S_i\subseteq T_i$. Since $G_1$ and $G_2$ are edge-disjoint, the pattern of this partition of $G$ is obtained by identifying one vertex of the pattern of the tree-partition of $G_1$ with one vertex of the pattern of the tree-partition of $G_2$. Since the patterns of the tree-partitions of $G_1$ and $G_2$ are forests, the pattern of the partition of $G$ is a forest, and we have a tree-partition of $G$. 

Moreover, $S$ is contained in a single bag $T_1\cup T_2$ and 
\begin{align*}
|T_1\cup T_2|
&=|T_1|+|T_2|-|Y|\\
&\leq\alpha|S_1|-\gamma(k+1)+\alpha|S_2|-\gamma(k+1)-(a+b)\\
&=\alpha(p_1+a+b)-\gamma(k+1)+\alpha(p_2+a+b)-\gamma(k+1)-(a+b)\\
&=\alpha(p_1+p_2+a)-2\gamma(k+1)+(\alpha-1)a+(2\alpha-1)b\\
&\leq\alpha|S|-2\gamma(k+1)+(2\alpha-1)(a+b)\\
&\leq\alpha|S|-2\gamma(k+1)+(2\alpha-1)(k+1)\\
&=\alpha|S|-\gamma(k+1)\enspace.
\end{align*}
Thus $|T_1\cup T_2|\leq
\alpha\cdot3(\gamma+1)(k+1)\Delta-\gamma(k+1)=
\gamma(k+1)(3\gamma\Delta-1)$. Since the other bags do not change we have the desired tree-partition of $G$.
\end{proof}

\section{General Lower Bound}
\seclabel{GeneralLowerBound}

The remainder of the paper studies lower bounds on the tree-partition-width. The graphs employed are chordal. We first show that tree-partitions of chordal graphs can be assumed to have certain useful properties. 


\begin{lemma}
\lemlabel{Chordal}
Every chordal graph $G$ has a tree-partition $T$ with width $\tpw{G}$, 
such that for every independent set $S$ of simplicial\footnote{A vertex is \emph{simplicial} if its neighbourhood is a clique.} vertices of $G$, 
and for every bag $B$ of $T$, 
either $B=\{v\}$ for some vertex $v\in S$, 
or the induced subgraph $G[B-S]$ is connected. 
\end{lemma}

\begin{proof}
Let $T_0$ be a tree-partition of a chordal graph $G$ with width $\tpw{G}$. 
Let $T$ be the partition of $G$ obtained from $T_0$ by replacing each bag $B$ of $T_0$ by bags corresponding to the connected components of $G[B]$. 
Then $T$ has width at most $\tpw{G}$. 

To prove that $T$ is a forest, 
suppose on the contrary that $T$ contains an induced cycle $C$. 
Since each bag in $C$ induces a connected subgraph of $G$, 
$G$ contains an induced cycle $D$ with at least one vertex from each bag in $C$. 
Since $G$ is chordal, $D$ is a triangle. 
Thus $C$ is a triangle, implying that the vertices in $D$ were in distinct bags in $T_0$ (since the bags of $T$ that replaced each bag of $T_0$ form an independent set).
Hence the bags of $T_0$ that contain $D$ induce a triangle in $T_0$, 
which is the desired contradiction since $T_0$ is a forest. 
Hence $T$ is a forest.

Let $S$ be an independent set of simplicial vertices of $G$.
Consider a bag $B$ of $T$. By construction, $G[B]$ is connected. 
First suppose that $B\subseteq S$. 
Since $S$ is an independent set and $G[B]$ is connected, 
$B=\{v\}$ for some vertex $v\in S$. 

Now assume that $B-S\ne\emptyset$. 
Suppose on the contrary that $G[B-S]$ is disconnected.
Thus $B\cap S$ is a cut-set in $G[B]$. Let $v$ and $w$ be vertices in distinct components of $G[B-S]$ such that the distance between $v$ and $w$ in $G[B]$ is minimised. (This is well-defined since $G[B]$ is connected.)\ 
Since $S$ is an independent set, every shortest path between $v$ and $w$ in $G[B]$ has only two edges. 
That is, $v$ and $w$ have a common neighbour $x$ in $B\cap S$. 
Since $x$ is simplicial, $v$ and $w$ are adjacent. 
This contradiction proves that $G[B-S]$ is connected. 
\end{proof}


The next lemma is the key component of the proof of \thmref{GeneralLowerBound}.
For integers $a<b$, let $[a,b]:=\{a,a+1,\dots,b\}$ and $[b]:=[1,b]$.

\begin{lemma}
\lemlabel{GeneralLowerBound}
For all integers $k\geq2$ and $\Delta\geq3k+1$, for infinitely many values of $N$ there is a chordal graph $G$ with $N$ vertices, tree-width $\tw{G}=2k-1$, maximum degree $\Delta(G)\leq\Delta$, and tree-partition-width $\tpw{G}>\quarter k(\Delta-3k)$.
\end{lemma}

\begin{proof}
Let $n$ be an integer with $n>\max\{\half k(\Delta-3k),2\}$. Let $H$ be the graph with vertex set $\{(x,y):x\in[n],y\in[k]\}$, where distinct vertices $(x_1,y_1)$ and $(x_2,y_2)$ are adjacent if and only if $|x_1-x_2|\leq 1$. The set of vertices $\{(x,y):y\in[k]\}$ is the \emph{$x$-column}. The set of vertices $\{(x,y):x\in[n]\}$ is the \emph{$y$-row}. Observe that each column induces a $k$-vertex clique, and each row induces an $n$-vertex path. 

Let $C$ be an induced cycle in $H$. If $(x,y)$ is a vertex in $C$ with $x$ minimum then the two neighbours of $(x,y)$ in $C$ are adjacent. Thus $C$ is a triangle. Hence $H$ is chordal. Observe that each pair of consecutive columns form a maximum clique of $2k$ vertices in $H$. Thus $H$ has tree-width $2k-1$.
Also note that $H$ has maximum degree $3k-1$. 


An edge of $H$ between vertices $(x,y)$ and $(x+1,y)$ is \emph{horizontal}. As illustrated in \figref{Grid}, construct a graph $G$ from $H$ as follows. For each horizontal edge $vw$ of $H$, add \ceil{\half(\Delta-3k)} new vertices, each adjacent to $v$ and $w$. Since $H$ is chordal and each new vertex is simplicial, $G$ is chordal. The addition of degree-$2$ vertices to $H$ does not increase the maximum clique size (since $k\geq2$). Thus $G$ has clique number $2k$ and tree-width $2k-1$.  Since each vertex of $H$ is incident to at most two horizontal edges, $G$ has maximum degree $3k-1+2\ceil{\half(\Delta-3k)}\leq\Delta$. 


\Figure{Grid}{\includegraphics{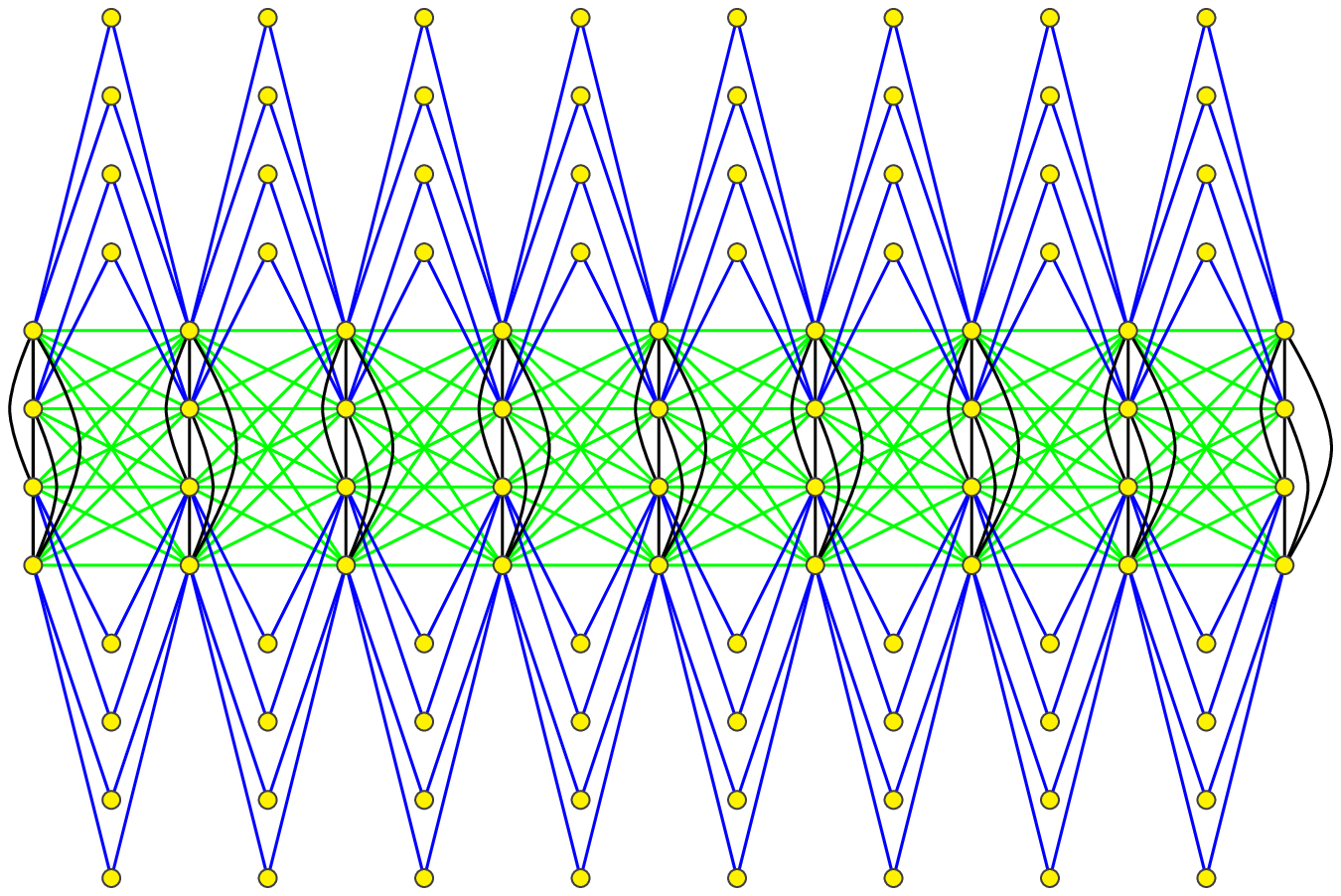}}{The graph $G$ with $k=4$, $\Delta=15$, and $n=9$.}

Observe that $V(G)-V(H)$ is an independent set of simplicial vertices in $G$. 
By \lemref{Chordal}, $G$ has a tree-partition $T$ with width $\tpw{G}$, 
such that for every bag $B$ of $T$, 
either $B=\{v\}$ for some vertex $v$ of $G-H$, 
or the induced subgraph $H[B]$ is connected. 
Since $G$ is connected, $T$ is a (connected) tree. 
Let $U$ be the tree-partition of $H$ induced by $T$. 
That is, to obtain $U$ from $T$ delete the vertices of $G-H$ from each bag, and delete empty bags. 
Since $H$ is connected, $U$ is a (connected) tree.
By \lemref{Chordal}, each bag of $U$ induces a connected subgraph of $H$.


Suppose that $U$ only has two bags $B$ and $C$. Then one of $B$ and $C$ contains at least $\half nk$ vertices. Since $k\geq2$, we have $\tpw{G}\geq\half nk>\quarter k(\Delta-3k)$, as desired. Now assume that $U$ has at least three bags. 

Consider a bag $B$ of $U$. Let $\ell(B)$ be the minimum integer such that some vertex in $B$ is in the $\ell(B)$-column, and let $r(B)$ be the maximum integer such that some vertex in $B$ is in the $r(B)$-column. Since $H[B]$ is connected, there is a path in $B$ from the $\ell(B)$-column to the $r(B)$-column. By the definition of $H$, for each $x\in[\ell(B),r(B)]$, the $x$-column contains a vertex in $B$. Let $I(B)$ be the closed real interval from $\ell(B)-\half$ to $r(B)+\half$. Observe that two bags $B$ and $C$ of $U$ are adjacent if and only if $I(B)\cap I(C)\neq\emptyset$. Thus $\{I(B):B\text{ is a bag of }U\}$ is an interval representation of the tree $U$. Every tree that is an interval graph is a caterpillar\footnote{A \emph{caterpillar} is a tree such that deleting the leaves gives a path.}; see \citep{Eckhoff-JGT93} for example. Thus $U$ is a caterpillar.

Let $\preceq$ be the relation on the set of non-leaf bags of $U$ defined by 
$A\preceq B$ if and only if $\ell(A)\leq\ell(B)$ and $r(A)\leq r(B)$. We claim that $\preceq$ is a total order. It is immediate that $\preceq$ is reflexive and transitive. To prove that $\preceq$ is antisymmetric, suppose on the contrary that $A\preceq B$ and $B\preceq A$ for distinct non-leaf bags $A$ and $B$. Thus $\ell(A)=\ell(B)$ and $r(A)=r(B)$. Since $U$ has at least three bags, there is a third bag $C$ that contains a vertex in the $(\ell(A)-1)$-column or in the $(r(A)+1)$-column. Thus $\{A,B,C\}$ induce a triangle in $U$, which is the desired contradiction. Hence $\preceq$ is antisymmetric. To prove that $\preceq$ is total, suppose on the contrary that $A\not\preceq B$ and $B\not\preceq A$ for distinct non-leaf bags $A$ and $B$. Now $A\not\preceq B$ implies that $\ell(A)>\ell(B)$ or $r(A)>r(B)$. Without loss of generality, $\ell(A)>\ell(B)$. Thus $B\not\preceq A$ implies that $r(B)>r(A)$. Hence the interval $[\ell(A),r(A)]$ is strictly within the interval $[\ell(B),r(B)]$ at both ends. For each $x\in[\ell(A),r(A)]$, every vertex in the $x$-column is in $A\cup B$, as otherwise $U$ would contain a triangle (since each column is a clique in $H$). Moreover, every vertex in the $(\ell(A)-1)$-column or in the $(r(A)+1)$-column is in $B$, as otherwise $U$ would contain a triangle (since the union of consecutive columns is a clique in $H$). Thus every neighbour of every vertex in $A$ is in $B$. That is, $A$ is a leaf in $U$. 
This contradiction proves that $\preceq$ is a total order on the set of non-leaf bags of $U$.

Suppose that $U$ has a $4$-vertex path $(A,B,C,D)$ as a subgraph. 

Thus $B$ and $C$ are non-leaf bags. Without loss of generality, $B\prec C$. If every column contains vertices in both $B$ and $C$, then $B$ and $C$ and any other bag would induce a triangle in $U$ (since each column induces a clique in $H$). Thus some column contains a vertex in $B$ but no vertex in $C$, and some column contains a vertex in $C$ but no vertex in $B$. Let $p$ be the maximum integer such that some vertex in $B$ is in the $p$-column, but no vertex in $C$ is in the $p$-column. Let $q$ be the minimum integer such that some vertex in $C$ is in the $q$-column, but no vertex in $B$ is in the $q$-column. Now $p<q$ since $B\prec C$. 

We claim that the $(p+1)$-column contains a vertex in $C$. If not, then the $(p+1)$-column contains no vertex in $B$ by the definition of $p$. Thus $r(B)=p$ since $H[B]$ is connected. Since $B$ is adjacent to $C$ in $U$, $\ell(C)\leq r(B)+1=p+1$. In particular, the $(p+1)$-column contains a vertex in $C$. Since $H[C]$ is connected, for $x\in[p+1,q]$, each $x$-column contains a vertex in $C$. In fact, $\ell(C)=p+1$ since the $p$-column contains no vertex in $C$. 
By symmetry, for $x\in[p,q-1]$, each $x$-column contains a vertex in $B$, and $r(C)=q-1$. 

The union of the $p$-column and the $(p+1)$-column only contains vertices in $B\cup C$, as otherwise $U$ would contain a triangle (since the union of two consecutive columns is a clique in $H$). By the definition of $p$, no vertex in the $p$-column is in $C$. Thus every vertex in the $p$-column is in $B$. By symmetry, every vertex in the $q$-column is in $C$. Now for each $y\in[k]$, the vertices $(p,y),(p+1,y),\dots,(q,y)$ are all in $B\cup C$, the first vertex $(p,y)$ is in $B$, and the last vertex $(q,y)$ is in $C$. Thus $(x,y)\in B$ and $(x+1,y)\in C$ for some $x\in[p,q-1]$. That is, in every row of $H$ there is a horizontal edge with one endpoint in $B$ and the other in $C$.

Thus there are at least $k$ horizontal edges with one endpoint in $B$ and the other in $C$ (now considered to be bags of $T$). For each such horizontal edge $vw$, each vertex of $G-H$ adjacent to $v$ and $w$ is in $B\cup C$, as otherwise $T$ would contain a triangle. There are \ceil{\half(\Delta-3k)} such vertices of $G-H$ for each of the $k$ horizontal edges between $B$ and $C$. Thus $|B\cup C|\geq\half k(\Delta-3k)$. Thus one of $B$ and $C$ has at least $\quarter k(\Delta-3k)$ vertices. Hence $\tpw{G}\geq\quarter k(\Delta-3k)$ as desired. 

Now assume that $U$ has no $4$-vertex path as a subgraph. 

A tree is a star if and only if it has no $4$-vertex path as a subgraph. Hence $U$ is a star. Let $R$ be the root bag of $U$. If $R$ contains a vertex in every column then $|R|\geq n$, implying $\tpw{G}\geq n\geq\quarter k(\Delta-3k)$, as desired. Now assume that for some $x\in[n]$, the $x$-column of $H$ contains no vertex in $R$. Let $B$ be a bag containing some vertex in the $x$-column. The $x$-column induces a clique in $H$, the only bag in $U$ that is adjacent to $B$ is $R$, and $R$ contains no vertex in the $x$-column. Thus every vertex in the $x$-column is in $B$. Since $R$ is the only bag in $U$ adjacent to $B$, there are at least $k$ horizontal edges with one endpoint in $B$ and the other endpoint in $R$. As in the case when $U$ contained a $4$-vertex path, we conclude that $\tpw{G}\geq\quarter k(\Delta-3k)$ as desired. 
\end{proof}

\begin{proof}[Proof of \thmref{GeneralLowerBound}]
Let $\ell:=\ceil{\frac{k}{2}}$. Thus $\ell\geq2$. By \lemref{GeneralLowerBound}, for each integer $\Delta\geq\Delta(k,\epsilon):=\max\{3\ell+1,\frac{3\ell}{8\epsilon}\}$, there are infinitely many values of $N$ for which there is a chordal graph $G$ with $N$ vertices, tree-width $\tw{G}=2\ell-1\leq k$, maximum degree $\Delta(G)\leq\Delta$, and tree-partition-width $\tpw{G}>\quarter \ell(\Delta-3\ell)$, which is at least $(\eighth-\epsilon)k\Delta$ since $\Delta\geq\frac{3\ell}{8\epsilon}$. 
\end{proof}

A \emph{domino} tree decomposition\footnote{See \citep{Diestel00} for an introduction to tree decompositions.} is a tree decomposition in which each vertex appears in at most two bags. The \emph{domino tree-width} of a graph $G$, denoted by \dtw{G}, is the minimum width of a domino tree decomposition of $G$.
Domino tree-width behaves like tree-partition-width in the sense that $\dtw{G}\geq\tw{G}$, and \dtw{G} is bounded for graphs of bounded tree-width and bounded degree \citep{BodEng-JAlg97}. The best upper bound is $$\dtw{G}\leq\big(9\,\tw{G}+7\big)\,\Delta(G)\,\big(\Delta(G)+1\big)-1,$$
which is due to \citet{Bodlaender-DMTCS99}, who also constructed a graph $G$ with $$\dtw{G}\geq\tfrac{1}{12}\,\tw{G}\,\Delta(G)-2.$$ 
Tree-partition-width and domino tree-width are related in that every graph $G$ satisfies $$\dtw{G}\geq\tpw{G}-1,$$as observed by \citet{BodEng-JAlg97}. Thus \thmref{GeneralLowerBound} provides examples of graphs $G$ with $$\dtw{G}\geq(\eighth-\epsilon)\,\tw{G}\,\Delta(G).$$ This represents a small constant-factor improvement over the above lower bound by \citet{Bodlaender-DMTCS99}.



\section{Lower Bound for Tree-width $2$}
\seclabel{SeriesParallel}

We now prove a lower bound on the tree-partition-width of graphs with tree-width $2$. 

\begin{theorem}
\thmlabel{LowerBound}
For all odd $\Delta\geq11$ there is a chordal graph $G$ with tree-width $2$,  maximum degree $\Delta$, and tree-partition-width $\tpw{G}\geq\twothirds(\Delta-1)$.
\end{theorem}

\begin{proof}
As illustrated in \figref{LowerBound}, 
let $G$ be the graph with $$V(G):=
\{r\}\cup\{v_i:i\in[\Delta]\}
\cup\{w_{i,\ell}:i\in[\Delta-1],\ell\in[\half(\Delta-3)]\}$$
and
$$E(G):=\{rv_i:i\in[\Delta]\}
\cup\{v_iv_{i+1}:i\in[\Delta-1]\}
\cup\{v_iw_{i,\ell},v_{i+1}w_{i,\ell}:i\in[\Delta-1],\ell\in[\half(\Delta-3)]\}.$$
Observe that $G$ has maximum degree $\Delta$. Clearly every induced cycle of $G$ is a triangle. Thus $G$ is chordal. Observe that $G$ has no $4$-vertex clique. Thus $G$ has tree-width $2$.

\begin{figure}[!h]
\begin{center}
\includegraphics{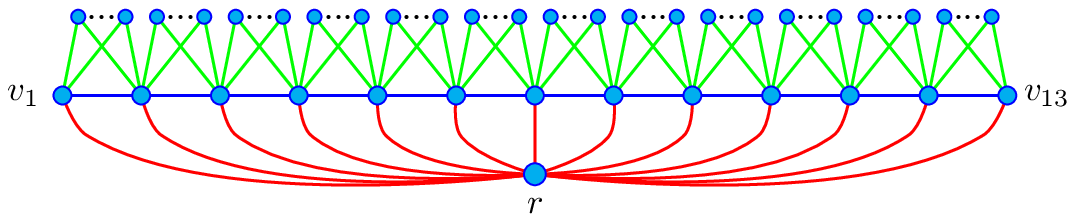}
\caption{\label{fig:LowerBound}Illustration for \thmref{LowerBound} with $\Delta=13$.}
\end{center}
\end{figure}

Let $T$ be the tree-partition of $G$ from \lemref{Chordal}. Then $T$ has width \tpw{G}, and every bag induces a connected subgraph of $G$. Let $R$ be the bag containing $r$. Let $B_1,\dots,B_d$ be the bags, not including $R$, that contain some vertex $v_i$. Thus $R$ is adjacent to each $B_j$ (since $r$ is adjacent to each $v_i$). Since $\{w_{i,\ell}:i\in[\Delta-1],\ell\in[\half(\Delta-3)]\}$ is an independent set of simplicial vertices, by \lemref{Chordal}, for each $j\in[d]$, the vertices $\{v_1,v_2,\dots,v_\Delta\}\cap B_j$ induce a (connected) subpath of $G$. 

First suppose that $d=0$. Then the $\Delta+1$ vertices $\{r,v_1,\dots,v_\Delta\}$ are contained in one bag $R$. Thus $\tpw{G}\geq\Delta+1\geq\twothirds(\Delta-1)$. 

Now suppose that $d=1$. Thus $\{r,v_1,\dots,v_\Delta\}\subseteq R\cup B_1$. In addition, at least one edge $v_iv_{i+1}$ has one endpoint in $R$ and the other endpoint in $B_1$. Thus $w_{i,\ell}\in R\cup B_1$ for each $\ell\in[\half(\Delta-3)\}]$. Hence $1+\Delta+\half(\Delta-3)$ vertices are contained in two bags. Thus one bag contains at least $\quarter(3\Delta-1)$ vertices, and $\tpw{G}\geq\quarter(3\Delta-1)\geq\twothirds(\Delta-1)$. 
 
Finally suppose that $d\geq2$. Since $\{v_1,v_2,\dots,v_\Delta\}\cap B_j$ induce a subpath in each bag $B_j$, we can assume that $\{v_1,v_2,\dots,v_\Delta\}\cap B_j=\{v_i:i\in[f(j),g(j)]\}$, where
$$1\leq 
f(1)\leq g(1)<
f(2)\leq g(2)<
\dots
<f(d)\leq g(d)\leq\Delta.$$

Distinct $B_j$ bags are not adjacent (since $T$ is a tree).
Thus $v_{f(j)-1}\in R$ for each $j\in[2,d]$. 
Similarly, $v_{g(j)+1}\in R$ for each $j\in[d-1]$.
Thus $w_{f(j)-1,\ell}\in R\cup B_j$ for each $j\in[2,d]$ and 
$\ell\in[\half(\Delta-3)\}]$. Similarly, $w_{g(j),\ell}\in R\cup B_j$ for each $j\in[d-1]$ and $\ell\in[\half(\Delta-3)\}]$. 

Hence the bags $R,B_1,\dots,B_d$ contain at least $$1+\Delta+2(d-1)\cdot\half(\Delta-3)$$ vertices. Therefore one of these bags has at least $$(1+\Delta+(d-1)(\Delta-3))/(d+1)$$ vertices, which is at least $\twothirds(\Delta-1)$. Hence $\tpw{G}\geq\twothirds(\Delta-1)$. 
\end{proof}

\begin{figure}[!h]
\begin{center}
\includegraphics{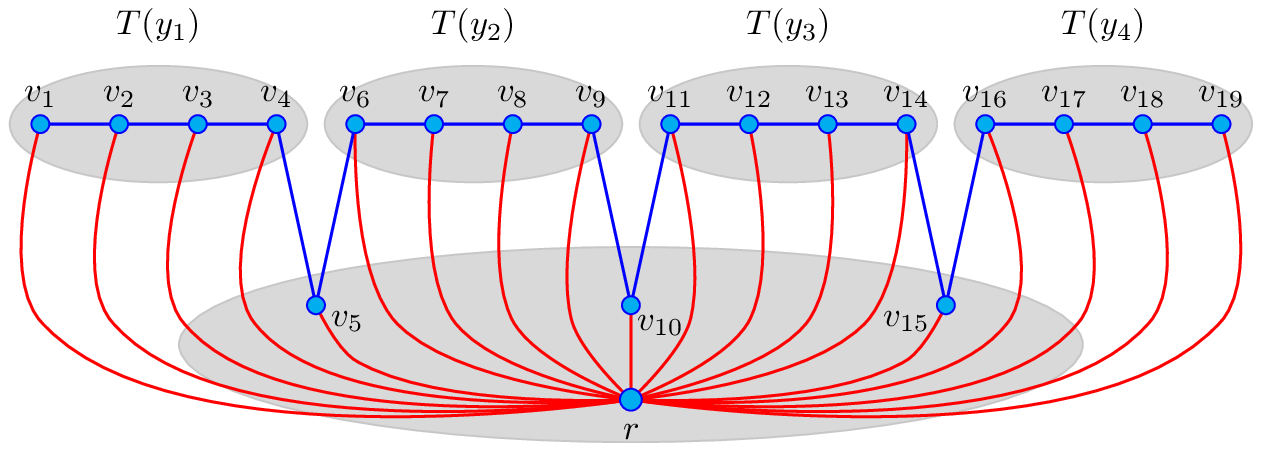}
\caption{Illustration for \thmref{LowerBound} with $\Delta=19$ and $d=4$.}
\end{center}
\end{figure}


\def\soft#1{\leavevmode\setbox0=\hbox{h}\dimen7=\ht0\advance \dimen7
  by-1ex\relax\if t#1\relax\rlap{\raise.6\dimen7
  \hbox{\kern.3ex\char'47}}#1\relax\else\if T#1\relax
  \rlap{\raise.5\dimen7\hbox{\kern1.3ex\char'47}}#1\relax \else\if
  d#1\relax\rlap{\raise.5\dimen7\hbox{\kern.9ex \char'47}}#1\relax\else\if
  D#1\relax\rlap{\raise.5\dimen7 \hbox{\kern1.4ex\char'47}}#1\relax\else\if
  l#1\relax \rlap{\raise.5\dimen7\hbox{\kern.4ex\char'47}}#1\relax \else\if
  L#1\relax\rlap{\raise.5\dimen7\hbox{\kern.7ex
  \char'47}}#1\relax\else\message{accent \string\soft \space #1 not
  defined!}#1\relax\fi\fi\fi\fi\fi\fi} \def\Dbar{\leavevmode\lower.6ex\hbox to
  0pt{\hskip-.23ex \accent"16\hss}D}

\end{document}